\documentclass[aps,pre,preprint,groupedaddress,floatfix]{revtex4-2}

\usepackage{bm,amsmath,amssymb,amsthm}
\usepackage[dvipdfm]{graphicx}
\usepackage{color}

\newtheorem{theorem}{Theorem}[section]
\newtheorem{lemma}[theorem]{Lemma}
\newtheorem{remark}{Remark}

\begin{document}


\title{Local H theorem for Enskog and Enskog--Vlasov equations with a modified Enskog factor}


\author{Aoto Takahashi}
\email[]{takahashi.aoto.63c@st.kyoto-u.ac.jp}
\author{Shigeru Takata}
\email[]{takata.shigeru.4a@kyoto-u.ac.jp}
\affiliation{Department of Aeronautics and Astronautics, Kyoto University,
Kyoto-daigaku-katsura, Kyoto 615-8540, Japan}


\date{\today}
\begin{abstract}

The local H theorem  is shown to hold for the Enskog equation with a modified Enskog factor proposed by the authors [Phys. Rev. E \textbf{111}, 065108 (2025)]. This is a stronger statement than the global one in the same paper and has been obtained along the lines of Mareschal \emph{et al.} [Phys. Rev. Lett. \textbf{52}, 1169--1172 (1984)] for the modified (or revised) Enskog equation. Furthermore, it is shown that the local H theorem also holds for the corresponding Enskog--Vlasov equation.

\end{abstract}

\maketitle
\section{Introduction}

The Enskog equation is a kinetic equation for dense gases. 
This equation devised by Enskog~\cite{E72} incorporates the center-of-mass displacement of colliding molecules into the Boltzmann equation and further includes the increased collision frequency resulting from the occupied volume of molecules.

The usefulness of the Enskog equation has become increasingly recognized with the development of recent numerical analyses, e.g., \cite{MS96,F97,MS97a,WLRZ16,F97b,HTT22,KAOMNFW25,HFNSV25}. 
However, the original equation by Enskog (OEE) retains the drawback that Boltzmann's H-theorem is not recovered due to the choice of the factor (the so-called Enskog factor) for increasing the collision frequency coming from the finite volume fraction effect. Fortunately, it was shown by Resibois~\cite{R78} that this drawback could be avoided for the modified (or revised) Enskog equation (MEE)~\cite{VE73}, but not for the OEE. 
Nevertheless, the mathematical intricacy has been hindering its further numerical applications so far.

In the meantime, we have recently proposed in \cite{TT25} a variant of the Enskog equation
that incorporates a slight modification to the Enskog factor that avoids the aforementioned drawback and allows numerical analysis with computational demands comparable to OEE. 
Indeed, we have shown that the H-theorem holds for EESM. 
However, the theorem there is a global statement for the total system quantity 
and is not a local statement that holds pointwisely in the physical system,
in contrast to the case of the Boltzmann equation.
Although the seminal work by Resibois~\cite{R78} was also a global statement for MEE, 
the local statement was later established by Mareschal{\it~et al.}~\cite{MBP84}.
In the present paper, along the lines of \cite{M83,MBP84}, 
we will show that the global statement of the H-theorem in \cite{TT25}
can be refined into a local statement.

The paper is organized as follows. 
First, the Enskog equation with a slight modification (EESM) 
is introduced in Sec.~\ref{sec:prbfrm}.
Then, in Sec.~\ref{sec:main} the main results for the local statement are presented. 
More specifically, using Appendix~\ref{app:shift} as an auxiliary, 
its kinetic part and collisional part are separately discussed
in Sec.~\ref{sec:kinH} and in Sec.~\ref{sec:colH},
and then these results are summarized to present the local theorem in Sec.~\ref{sec:LocalHtheorem}. 
Next, in Sec.~\ref{sec:EV}, the results of Sec.~\ref{sec:main} are extended to the Enskog--Vlasov equation that incorporates intermolecular attractive effects into EESM (EVESM, for short). 
The recovery of the global statement in \cite{TT25} from the present results
is addressed in Appendix~\ref{app:mono}. 
Finally, the paper is concluded in Sec.~\ref{sec:Conclusion}.

\section{Enskog equation with a slight modification (EESM) \label{sec:prbfrm}}

We consider the Enskog equation for a single species dense gas
that is composed of hard sphere molecules with a common diameter $\sigma$ and mass $m$.
Let $D$ be a fixed spatial domain in which the center of gas molecules is confined. 
Let $t$, $\bm{X}$, and $\bm{\xi}$
be a time, a spatial position, and a molecular velocity, respectively.
Denoting the one-particle distribution function of gas molecules
by $f(t,\bm{X},\bm{\xi}$), 
the Enskog equation is written as
\begin{subequations}\label{MEE}
\begin{align}
 & \frac{\partial f}{\partial t}+\xi_{i}\frac{\partial f}{\partial X_{i}}=J(f)\equiv J^{G}(f)-J^{L}(f),\quad \mathrm{for\ }\bm{X}\in D,\displaybreak[0]\label{eq:2.1}\\
 & J^{G}(f)\equiv\frac{\sigma^{2}}{m}\int {g(\bm{X}_{\sigma\bm{\alpha}}^{+},\bm{X})f_{*}^{\prime}(\bm{X}_{\sigma\bm{\alpha}}^{+})f^{\prime}(\bm{X})} \notag\\&\qquad\qquad\times V_{\alpha}\theta(V_{\alpha})d\Omega(\bm{\alpha})d\bm{\xi}_{*},\displaybreak[0]\label{eq:2.2}\\
 & J^{L}(f)\equiv\frac{\sigma^{2}}{m}\int {g(\bm{X}_{\sigma\bm{\alpha}}^{-},\bm{X})f_{*}(\bm{X}_{\sigma\bm{\alpha}}^{-})f(\bm{X})} \notag\\&\qquad\qquad\times V_{\alpha}\theta(V_{\alpha})d\Omega(\bm{\alpha})d\bm{\xi}_{*},\label{eq:2.3}
\end{align}
\end{subequations}
\noindent
where $\bm{X}_{\bm{x}}^{\pm}=\bm{X}\pm\bm{x}$,
$\bm{\alpha}$ is a unit vector, 
$d\Omega(\bm{\alpha})$ is a
solid angle element in the direction of $\bm{\alpha}$,
$\theta$ is the Heaviside function
\begin{subequations}
\begin{equation}
\theta(x)=\begin{cases}
1, & x\ge0\\
0, & x<0
\end{cases},
\end{equation}
\noindent
and the following notation convention has been used:
\begin{align}
 & \begin{cases}
{
f(\bm{X})=f(\bm{X},\bm{\xi}),\ f^{\prime}(\bm{X})=f(\bm{X},\bm{\xi}^{\prime})},\\
{
f_{*}(\bm{X}_{\sigma\bm{\alpha}}^{-})=f(\bm{X}_{\sigma\bm{\alpha}}^{-},\bm{\xi}_{*}),\ f_{*}^{\prime}(\bm{X}_{\sigma\bm{\alpha}}^{+})=f(\bm{X}_{\sigma\bm{\alpha}}^{+},\bm{\xi}_{*}^{\prime})},
\end{cases}\label{eq:contf}\displaybreak[0]\\
 & \begin{cases}
{
 \bm{\xi}^{\prime}=\bm{\xi}+V_{\alpha}\bm{\alpha},\quad
 \bm{\xi}_{*}^{\prime}=\bm{\xi}_{*}-V_{\alpha}\bm{\alpha},
 }\\
{ 
 V_{\alpha}=\bm{V}\cdot\bm{\alpha},\quad\bm{V}=\bm{\xi_{*}}-\bm{\xi}.}
 \end{cases}\label{eq:2.5}
\end{align}
\noindent
Here and in what follows, the argument $t$ is often suppressed, unless confusion is anticipated. 
The convention \eqref{eq:contf} will be applied only to the quantities that depend on molecular velocity.
It should be noted that \eqref{eq:2.2} [or \eqref{eq:2.3}] makes sense only 
when $\bm{X}^+_{\sigma\bm{\alpha}}$
[or $\bm{X}^-_{\sigma\bm{\alpha}}$] as well as $\bm{X}$ is in the domain $D$.
Accordingly, the factor $g$ should be understood as
\begin{equation}
g(\bm{X},\bm{Y})=\mathsf{g}(\bm{X},\bm{Y})\chi_D(\bm{X})\chi_D(\bm{Y}),\label{eq:sfg}
\end{equation}%
%
\noindent
where $\chi_D$ is the indicator function:
\begin{equation}
\chi_{D}(\bm{X})=\begin{cases}
1, & \bm{X}\in D\\
0, & \mbox{otherwise}
\end{cases},\label{eq:chi_def}
\end{equation}
\end{subequations}
\noindent
so that the integration in
(\ref{eq:2.2}) and (\ref{eq:2.3}) is over the entire space
of $\bm{\xi}_*$ and for all directions of $\bm{\alpha}$, irrespective of the position in the domain $D$.
The range of integration with respect to $\bm{\alpha}$, $\bm{\xi}$, and $\bm{\xi}_*$
will be suppressed in most cases, unless confusion is anticipated in this paper.

The factor $\mathsf{g}$ occurring in \eqref{eq:sfg} 
is the so-called Enskog factor and is generically
assumed to be positive and symmetric with respect to the exchange of two position vectors: $\mathsf{g}(\bm{X},\bm{Y})=\mathsf{g}(\bm{Y},\bm{X})$.
Although there are some varieties of $\mathsf{g}$ in the literature, e.g., \cite{E72,VE73,DVK21,BLPT91,BB18,BB19,TT25},
our target is the following one proposed in \cite{TT25}:
\begin{subequations}\label{eq:g_prop}\begin{align}
& \mathsf{g}(\bm{X},\bm{Y})=\mathcal{S}(\mathcal{R}(\bm{X}))+\mathcal{S}(\mathcal{R}(\bm{Y})),
\label{eq:g_def}
\\
& \mathcal{R}(\bm{X})=\frac{1}{m}\int_D \rho(\bm{Y})\theta(\sigma-|\bm{Y}-\bm{X}|) d\bm{Y},
\label{eq:R_def}
\\
& \rho(\bm{X})=\langle f(\bm{X}) \rangle,
\quad \langle \bullet \rangle=\int_{\mathbb{R}^3} \bullet \, d\bm{\xi},
\label{eq:rho_def}
\end{align}\end{subequations}
\noindent
where $\mathcal{S}$ is a non-negative function, 
the form of which is determined according to the target gas 
as explained in the next paragraph.
The Enskog equation incorporating \eqref{eq:g_prop}
is referred to as EESM in \cite{TT25}, and this terminology is also used in the present paper. 

As explained in \cite{TT25},
the equation of state (EoS) at the uniform equilibrium state
is expressed as
\begin{equation}
p=\rho RT (1+2b\rho \mathcal{S}(2b\rho)),
\quad
b=\frac{2\pi}{3}\frac{\sigma^3}{m},
\end{equation}
\noindent
where $p$ is the pressure of the gas.
Hence, for instance, by setting
\begin{subequations}%
\begin{equation}
\mathcal{S}(x)=\frac{1}{2-x}, \label{eq:SvdW}
\end{equation}
\noindent
the van der Waals EoS \cite{vdW} for non-attractive molecules 
\begin{equation}
p=\frac{\rho RT}{1-b\rho}=\rho RT(1+\frac{b\rho}{1-b\rho}),
\label{eq:EoSvdW}
\end{equation}\end{subequations}
\noindent
is recovered,
while by setting
\begin{subequations}%
\begin{equation}
 \mathcal{S}(x)=\frac{16(16-x)}{(8-x)^3}, \label{eq:SCS}
\end{equation}
\noindent
the Carnahan--Starling EoS \cite{CS69}
\begin{equation}
p=\rho RT\frac{1+\eta+\eta^2-\eta^3}{(1-\eta)^3}
 =\rho RT(1+\frac{4\eta-2\eta^2}{(1-\eta)^3}),
\label{eq:EoSCS}
\end{equation}\end{subequations}
\noindent
is recovered,
where $\eta={b\rho}/{4}$.

It should be noted that
the collision term of the Enskog equation is not responsible 
for the attractive part of the {EoS}.
The attractive part is to be recovered by the Vlasov term of the Enskog--Vlasov equation; 
see, e.g., \cite{G71,FGLS18}.
Therefore the form \eqref{eq:SvdW} applies to the {EESM}
with the Vlasov term for the \emph{full} version of the van der Waals fluids as well, see Sec.~\ref{sec:VlasEoS}.

In closing this section, let us list the definitions
of macroscopic quantities for later convenience and summarize the conservation laws.
In addition to the density $\rho$ already given in \eqref{eq:rho_def},
the flow velocity $\bm{v}$ (or $v_i$) and temperature $T$
are defined by
\begin{subequations}
\begin{align}
& v_i=\frac{1}{\rho}\langle \xi_i f \rangle, \label{eq:flow_def}\\
& T=\frac{1}{3R\rho}\langle (\bm{\xi}-\bm{v})^2 f \rangle, \label{eq:temp_def}
\end{align}
\noindent
and the so-called kinetic part of the specific internal energy $e^{(k)}$,
that of the stress tensor $p_{ij}^{(k)}$, and that of the heat-flow vector $\bm{q}^{(k)}$ (or $q_i^{(k)}$), are defined by
\begin{align}
& e^{(k)}=\frac{1}{2\rho}\langle (\bm{\xi}-\bm{v})^2 f \rangle
  (=\frac32RT), \label{eq:kin_inteng_def}\\
& p_{ij}^{(k)}=\langle(\xi_i-v_i)(\xi_j-v_j)f\rangle, \label{eq:kin_stress_def}\\
& q_i^{(k)}=\frac12\langle(\xi_i-v_i)(\bm{\xi}-\bm{v})^2 f\rangle. \label{eq:kin_heatfl_def}
\end{align}
\end{subequations}
\noindent
The conservation laws that are derived from the Enskog equation are 
\begin{subequations}
\begin{align}
&
\frac{\partial\rho}{\partial t}+\frac{\partial\rho v_i}{\partial X_i}=0,\label{eq:continuity} \\
&
\frac{\partial\rho v_j}{\partial t}
+\frac{\partial}{\partial X_i}(\rho v_i v_j+p_{ij})
=0, \label{eq:mom_div} \\
&
\frac{\partial}{\partial t}[\rho(e+\frac12 \bm{v}^2)]
+\frac{\partial}{\partial X_i}[\rho v_i(e+\frac12 \bm{v}^2)
+p_{ij}v_j+q_i]=0,
\label{eq:eng_div}
\end{align}\label{eq:conservation}
\end{subequations}
\noindent
where
\begin{subequations}
\begin{equation}
p_{ij}=p_{ij}^{(k)}+p_{ij}^{(c)},\quad
q_i=q_i^{(k)}+q_i^{(c)},\quad
e=e^{(k)},\label{eq:kin+col}
\end{equation}
\begin{align}
p_{ij}^{(c)}= & \frac{\sigma^{2}}{2m}\int {\int_{0}^{\sigma}}\alpha_{i}\alpha_{j}V_{\alpha}^{2}\theta(V_{\alpha}) g(\bm{X}_{\lambda\bm{\alpha}}^{+},\bm{X}_{(\lambda-\sigma)\bm{\alpha}}^{+}) \notag\\&\quad\times f_{*}(\bm{X}_{(\lambda-\sigma)\bm{\alpha}}^{+})f(\bm{X}_{\lambda\bm{\alpha}}^{+}){d\lambda} d\Omega(\bm{\alpha})d\bm{\xi}_{*}d\bm{\xi}, \label{eq:col_stress}\\
q_{i}^{(c)}= & \frac{\sigma^{2}}{4m}\int{\int_{0}^{\sigma}}\alpha_{i}[(\bm{c}+\bm{c}_{*})\cdot\bm{\alpha}]V_{\alpha}^{2}\theta(V_{\alpha}) \notag\\&\quad\times g(\bm{X}_{\lambda\bm{\alpha}}^{+},\bm{X}_{(\lambda-\sigma)\bm{\alpha}}^{+})f_{*}(\bm{X}_{(\lambda-\sigma)\bm{\alpha}}^{+}) \notag\\
&\qquad\times f(\bm{X}_{\lambda\bm{\alpha}}^{+}){d\lambda} d\Omega(\bm{\alpha})d\bm{\xi}_{*}d\bm{\xi},\label{eq:4.6}
\end{align}\label{eq:constitution}%
\end{subequations}
\noindent
$\bm{c}=\bm{\xi}-\bm{v}$, and $\bm{c}_{*}=\bm{\xi}_{*}-\bm{v}$;
see, e.g., \cite{CL88,F99}.
$p_{ij}^{(c)}$ and $\bm{q}^{(c)}$ (or $q_i^{(c)}$) are the so-called collisional part of
the stress tensor and heat-flow vector, respectively.
The reader is referred to \cite{TT25} for the derivation of \eqref{eq:conservation} and \eqref{eq:constitution}.
  
\section{Main results\label{sec:main}}
\subsection{Kinetic part of the local H function}\label{sec:kinH}

First we shall focus on the so-called kinetic part of the local H function
that is defined by
\begin{equation}
{H}^{(k)}\equiv\langle f{\ln f} \rangle.\label{H_kinetic}
\end{equation}
\noindent
Multiplying the Enskog equation (\ref{eq:2.1}) by $(1+{\ln f})$
and integrating the result with respect to $\bm{\xi}$ gives
\begin{equation}
 \frac{\partial}{\partial t}\langle f{\ln f}\rangle
+\frac{\partial}{\partial X_{i}}\langle\xi_{i}f{\ln f}\rangle
=\langle J(f){\ln f}\rangle.\label{eq:27}
\end{equation}

As is explained in Appendix~\ref{app:shift},
the right-hand side of \eqref{eq:27} can be transformed into the form that
\begin{align}
 & \langle J(f){\ln f}\rangle\nonumber \\
= & \frac{\sigma^{2}}{2m}\int\ln\Big(\frac{f_{*}^{\prime}(\bm{X}_{\sigma\bm{\alpha}}^{-})f^{\prime}(\bm{X})}{f_{*}(\bm{X}_{\sigma\bm{\alpha}}^{-})f(\bm{X})}\Big)
g(\bm{X},\bm{X}_{\sigma\bm{\alpha}}^{-}) \notag \\
&\times f(\bm{X})f_{*}(\bm{X}_{\sigma\bm{\alpha}}^{-})V_{\alpha}\theta(V_{\alpha})d\Omega(\bm{\alpha})d\bm{\xi}d\bm{\xi}_{*}
-\frac{\partial J_i^{(k)}}{\partial X_i}.
\label{eq:Jlnf}
\end{align}
\noindent
Note that the second term on the right-hand side was absent in \cite{TT25}, since
this term vanishes by integration in space.
Since $x\ln(y/x)\le y-x$ for any $x,y>0$, 
we obtain from \eqref{eq:Jlnf} that
\begin{subequations}%
\begin{equation}
\langle J(f){\ln f}\rangle \le \mathcal{I}-\frac{\partial J_i^{(k)}}{\partial X_i},
\label{eq:Jlnf_1}
\end{equation}
\noindent
where 
\begin{align}
\mathcal{I}&=\frac{\sigma^{2}}{2m}\int g(\bm{X},\bm{X}_{\sigma\bm{\alpha}}^{-})[f_{*}^{\prime}(\bm{X}_{\sigma\bm{\alpha}}^{-})f^{\prime}(\bm{X}) \notag\\&\quad -f(\bm{X})f_{*}(\bm{X}_{\sigma\bm{\alpha}}^{-})] V_{\alpha}\theta(V_{\alpha})d\Omega(\bm{\alpha})d\bm{\xi}d\bm{\xi}_{*},\label{eq:I_def}
\end{align}
\noindent
and the equality in \eqref{eq:Jlnf_1} holds if and only if 
\begin{equation}
 f_{*}^{\prime}(\bm{X}_{\sigma\bm{\alpha}}^{-})f^{\prime}(\bm{X})
-f(\bm{X})f_{*}(\bm{X}_{\sigma\bm{\alpha}}^{-})=0, \label{eq:equilibrium}
\end{equation}%
\end{subequations}
\noindent
or equivalently $\mathcal{I}=0$. 
Therefore, it holds that
\begin{equation}
 \frac{\partial}{\partial t} \langle f\ln f\rangle
+\frac{\partial}{\partial X_i} \Big( \langle\xi_i f\ln f \rangle +J_i^{(k)}\Big) 
\le \mathcal{I}. 
\label{eq:Hk_inequality}
\end{equation}
As is explained in Appendix~\ref{app:ColEnt},
$\mathcal{I}$ is eventually reduced to a much simpler form:
\begin{align}
\mathcal{I}=
&\frac{\sigma^{2}}{2m}\int g(\bm{X},\bm{X}_{\sigma\bm{\alpha}}^{+})
 \rho(\bm{X})\rho(\bm{X}_{\sigma\bm{\alpha}}^{+})\notag\\
&\qquad
\times[ v_j(\bm{X}_{\sigma\bm{\alpha}}^{+})-v_j(\bm{X})]\alpha_j 
 d\Omega(\bm{\alpha}).
\label{eq:I_final}
\end{align}
\noindent
This simpler form will be useful in Sec.~\ref{sec:colH}. 

\begin{remark}
Since the condition \eqref{eq:equilibrium} is equivalent to
\begin{equation}
 \ln f         (\bm{X})
+\ln f_*       (\bm{X}^-_{\sigma\bm{\alpha}})
=\ln f^\prime  (\bm{X})
+\ln f_*^\prime(\bm{X}^-_{\sigma\bm{\alpha}}),
\end{equation}
\noindent
the $\ln f$ satisfying this equation 
is the summational invariant. Hence, $f$ is restricted to the form
\begin{subequations}\label{eq:Max}
\begin{align}
 f(\bm{X})
=&\frac{\rho(\bm{X})}{(2\pi RT)^{3/2}}
\exp(-\frac{[\bm{\xi}-\bm{v}(\bm{X})]^2}{2RT}), \\
\bm{v}(\bm{X})=&\bm{u}+\bm{X}\times\bm{\omega},
\end{align}
\end{subequations}
\noindent
see, e.g., \cite{TT24,MGB18}.
Note that $T$, $\bm{u}$, and $\bm{\omega}$ in \eqref{eq:Max} 
are independent of $\bm{X}$ and that
$\bm{u}$ and $\bm{\omega}$ represent
the translational and the angular velocity of the flow, respectively.
\end{remark}

\subsection{Collisional part of the local H function\label{sec:colH}}

Next consider the function
\begin{equation}
{H}^{(c)}= \rho(\bm{X}) [\int_0^{\mathcal{R}(\bm{X})}
\mathcal{S}(x)dx ]. \label{eq:H^c}
\end{equation}
\noindent
Using a concise notation $\bm{r}=\bm{Y}-\bm{X}$ and $r=|\bm{r}|$,
its time derivative is transformed as
\begin{widetext}%
\begin{align}
\frac{\partial}{\partial t}{H}^{(c)}
=& \frac{\partial\rho(\bm{X})}{\partial t}
 \int_0^{\mathcal{R}(\bm{X})}\mathcal{S}(x)dx 
+\rho(\bm{X})\frac{\partial\mathcal{R}(\bm{X})}{\partial t}
\mathcal{S}(\mathcal{R}(\bm{X})) \notag \displaybreak[0]\\
=&\frac{\partial\rho(\bm{X})}{\partial t}
  \int_0^{\mathcal{R}(\bm{X})}\mathcal{S}(x)dx
+ \frac{\rho(\bm{X})}{m}[\int_D\frac{\partial\rho(\bm{Y})}{\partial t}
  \theta(\sigma-r)d\bm{Y}]
  \mathcal{S}(\mathcal{R}(\bm{X})) \notag \displaybreak[0]\\
=&-\frac{\partial\rho(\bm{X}) v_i(\bm{X})}{\partial X_i}
  \int_0^{\mathcal{R}(\bm{X})}\mathcal{S}(x)dx
 -\frac{\rho(\bm{X})}{m}[\int_D\frac{\partial\rho(\bm{Y}) v_i(\bm{Y})}{\partial Y_i}
  \theta(\sigma-r)d\bm{Y}]
  \mathcal{S}(\mathcal{R}(\bm{X})) \notag \displaybreak[0]\\
=&-\frac{\partial}{\partial X_i}
   \Big( \rho(\bm{X}) v_i(\bm{X})\int_0^{\mathcal{R}(\bm{X})}\mathcal{S}(x)dx \Big)
  +\rho(\bm{X}) v_i(\bm{X})
   \frac{\partial\mathcal{R}(\bm{X})}{\partial X_i}\mathcal{S}(\mathcal{R}(\bm{X}))\notag \\
 &-\frac{\rho(\bm{X})}{m}[\int_D\frac{\partial}{\partial Y_i}
   \Big(\rho(\bm{Y}) v_i(\bm{Y})\theta(\sigma-r)\Big)d\bm{Y}] 
   \mathcal{S}(\mathcal{R}(\bm{X})) \notag\\
 &+\frac{\rho(\bm{X})}{m}[\int_D\rho(\bm{Y}) v_i(\bm{Y})\frac{\partial}{\partial Y_i}
   \theta(\sigma-r)d\bm{Y}]
   \mathcal{S}(\mathcal{R}(\bm{X})) \notag \displaybreak[0]\\
=&-\frac{\partial}{\partial X_i}
   \Big( \rho(\bm{X}) v_i(\bm{X})\int_0^{\mathcal{R}(\bm{X})}\mathcal{S}(x)dx \Big)
  +\rho(\bm{X}) v_i(\bm{X})
   \int_D\frac{\rho(\bm{Y})}{m}\delta(\sigma-r)\frac{r_i}{r}d\bm{Y}
   \mathcal{S}(\mathcal{R}(\bm{X}))\notag \\
 &-\frac{\rho(\bm{X})}{m}[\int_D\frac{\partial}{\partial Y_i}
   \Big(\rho(\bm{Y}) v_i(\bm{Y})\theta(\sigma-r)\Big)d\bm{Y}] 
   \mathcal{S}(\mathcal{R}(\bm{X})) \notag\\
 &-\frac{\rho(\bm{X})}{m}[\int_D\rho(\bm{Y}) v_i(\bm{Y})\delta(\sigma-r)\frac{r_i}{r}d\bm{Y}]
   \mathcal{S}(\mathcal{R}(\bm{X})) \notag \displaybreak[0]\\
=&-\frac{\partial}{\partial X_i}
   \Big( H^{(c)}(\bm{X})v_i(\bm{X}) \Big)
-\frac{\rho(\bm{X})}{m}[\int_D\frac{\partial}{\partial Y_i}
   \Big(\rho(\bm{Y}) v_i(\bm{Y})\theta(\sigma-r)\Big)d\bm{Y}] 
   \mathcal{S}(\mathcal{R}(\bm{X})) 
\notag \\
 &+\rho(\bm{X}) v_i(\bm{X})
   \int_{\mathbb{R}^3}\frac{\rho(\bm{X}^+_{\bm{r}})}{m}\delta(\sigma-r)\frac{r_i}{r}
   \chi_D(\bm{X}^+_{\bm{r}})d\bm{r}
   \mathcal{S}(\mathcal{R}(\bm{X}))\notag\\
 &-\frac{\rho(\bm{X})}{m}[\int_{\mathbb{R}^3}\rho(\bm{X}^+_{\bm{r}}) v_i(\bm{X}^+_{\bm{r}})\delta(\sigma-r)\frac{r_i}{r}
   \chi_D(\bm{X}^+_{\bm{r}}) d\bm{r}]
   \mathcal{S}(\mathcal{R}(\bm{X})) \notag \displaybreak[0]\\
=&-\frac{\partial}{\partial X_i}
   \Big( H^{(c)}(\bm{X})v_i(\bm{X}) \Big)
-\frac{\rho(\bm{X})}{m}[\int_D\frac{\partial}{\partial Y_i}
   \Big(\rho(\bm{Y}) v_i(\bm{Y})\theta(\sigma-r)\Big)d\bm{Y}] 
   \mathcal{S}(\mathcal{R}(\bm{X})) \notag\\
 &-\frac{\sigma^2}{m}\int_{\mathbb{S}^2}
   [\rho(\bm{X})\rho(\bm{X}^+_{\sigma\bm{\alpha}}) v_i(\bm{X}^+_{\sigma\bm{\alpha}})
   -\rho(\bm{X}) v_i(\bm{X})\rho(\bm{X}^+_{\sigma\bm{\alpha}})]
    \chi_D(\bm{X}^+_{\sigma\bm{\alpha}})\alpha_i d\Omega(\bm{\alpha})
    \mathcal{S}(\mathcal{R}(\bm{X})).
\end{align}
\noindent
Since the last term is further transformed as
\begin{align}
 &-\frac{\sigma^2}{m}\int_{\mathbb{S}^2}
   [\rho(\bm{X})\rho(\bm{X}^+_{\sigma\bm{\alpha}}) v_i(\bm{X}^+_{\sigma\bm{\alpha}})
   -\rho(\bm{X}) v_i(\bm{X})\rho(\bm{X}^+_{\sigma\bm{\alpha}})]
   \chi_D(\bm{X})\chi_D(\bm{X}^+_{\sigma\bm{\alpha}})\alpha_id\Omega(\bm{\alpha})
   \mathcal{S}(\mathcal{R}(\bm{X})) \notag \displaybreak[0]\\
=&-\frac{\sigma^2}{m}\int_{\mathbb{S}^2}
   [v_i(\bm{X}^+_{\sigma\bm{\alpha}})
   -v_i(\bm{X})]\rho(\bm{X}^+_{\sigma\bm{\alpha}}) 
   \rho(\bm{X})\mathcal{S}(\mathcal{R}(\bm{X}))
   \chi_D(\bm{X})\chi_D(\bm{X}^+_{\sigma\bm{\alpha}})\alpha_i
   d\Omega(\bm{\alpha}) \notag \displaybreak[0]\\
=&-\frac{\sigma^2}{2m}\int_{\mathbb{S}^2}
   [v_i(\bm{X}^+_{\sigma\bm{\alpha}})
   -v_i(\bm{X})]\rho(\bm{X}^+_{\sigma\bm{\alpha}}) 
   \rho(\bm{X})[\mathcal{S}(\mathcal{R}(\bm{X}))+\mathcal{S}(\mathcal{R}(\bm{X}^+_{\sigma\bm{\alpha}}))]
   \chi_D(\bm{X})\chi_D(\bm{X}^+_{\sigma\bm{\alpha}})\alpha_i
   d\Omega(\bm{\alpha}) \notag \displaybreak[0]\\
&-\frac{\sigma^2}{2m}\int_{\mathbb{S}^2}
   [v_i(\bm{X}^+_{\sigma\bm{\alpha}})
   -v_i(\bm{X})]\rho(\bm{X}^+_{\sigma\bm{\alpha}}) 
   \rho(\bm{X})[\mathcal{S}(\mathcal{R}(\bm{X}))-\mathcal{S}(\mathcal{R}(\bm{X}^+_{\sigma\bm{\alpha}}))]
   \chi_D(\bm{X})\chi_D(\bm{X}^+_{\sigma\bm{\alpha}})\alpha_i
   d\Omega(\bm{\alpha}) \notag \displaybreak[0]\\
=&-\frac{\sigma^2}{2m}\int_{\mathbb{S}^2}
   [v_i(\bm{X}^+_{\sigma\bm{\alpha}})
   -v_i(\bm{X})]\rho(\bm{X}^+_{\sigma\bm{\alpha}}) 
   \rho(\bm{X})g(\bm{X},\bm{X}^+_{\sigma\bm{\alpha}})
   \alpha_id\Omega(\bm{\alpha}) \notag \displaybreak[0]\\
&+\frac{\sigma^2}{2m}\int_{\mathbb{S}^2}
   \Big(
   [v_i(\bm{X}^+_{\sigma\bm{\alpha}})
   -v_i(\bm{X})]\rho(\bm{X}^+_{\sigma\bm{\alpha}}) 
   \rho(\bm{X})\mathcal{S}(\mathcal{R}(\bm{X}^+_{\sigma\bm{\alpha}}))
   \chi_D(\bm{X})\chi_D(\bm{X}^+_{\sigma\bm{\alpha}})
   \notag \displaybreak[0]\\
& -[v_i(\bm{X})-v_i(\bm{X}^-_{\sigma\bm{\alpha}})]
   \rho(\bm{X}^-_{\sigma\bm{\alpha}}) 
   \rho(\bm{X})\mathcal{S}(\mathcal{R}(\bm{X}))
   \chi_D(\bm{X})\chi_D(\bm{X}^-_{\sigma\bm{\alpha}})
   \Big)\alpha_i
   d\Omega(\bm{\alpha}) \notag \displaybreak[0]\\
=&-\mathcal{I} 
  +\frac{\sigma^2}{2m}\int_{\mathbb{S}^2}\int_0^1
   \frac{\partial}{\partial\tau}\Big(
   [v_i(\bm{X}^+_{\tau\sigma\bm{\alpha}})
   -v_i(\bm{X}^+_{(\tau-1)\sigma\bm{\alpha}})]\rho(\bm{X}^+_{\tau\sigma\bm{\alpha}}) 
   \rho(\bm{X}^+_{(\tau-1)\sigma\bm{\alpha}})\notag\\
&\qquad\times  \mathcal{S}(\mathcal{R}(\bm{X}^+_{\tau\sigma\bm{\alpha}}))
   \chi_D(\bm{X}^+_{(\tau-1)\sigma\bm{\alpha}})\chi_D(\bm{X}^+_{\tau\sigma\bm{\alpha}})
   \Big)d\tau
   \alpha_i d\Omega(\bm{\alpha}) \notag \displaybreak[0]\\
=&-\mathcal{I} 
  +\frac{\partial}{\partial X_j}
   \Big(
   \frac{\sigma^3}{2m}\int_{\mathbb{S}^2}\int_0^1
   [v_i(\bm{X}^+_{\tau\sigma\bm{\alpha}})
   -v_i(\bm{X}^+_{(\tau-1)\sigma\bm{\alpha}})]\rho(\bm{X}^+_{\tau\sigma\bm{\alpha}}) 
   \rho(\bm{X}^+_{(\tau-1)\sigma\bm{\alpha}})\notag\\
&\qquad\times  \mathcal{S}(\mathcal{R}(\bm{X}^+_{\tau\sigma\bm{\alpha}}))
   \chi_D(\bm{X}^+_{(\tau-1)\sigma\bm{\alpha}})\chi_D(\bm{X}^+_{\tau\sigma\bm{\alpha}})
   d\tau
   \alpha_i \alpha_j d\Omega(\bm{\alpha}) \Big),\label{eq:Jcori}
\end{align}
\noindent
it holds that
\begin{subequations}%
\begin{equation}
\frac{\partial {H}^{(c)}}{\partial t}+\frac{\partial J_i^{(c)}}{\partial X_i}
=-\mathcal{I}
 -\frac{\rho(\bm{X})}{m}[\int_D\frac{\partial}{\partial Y_i}
   \Big(\rho(\bm{Y}) v_i(\bm{Y})\theta(\sigma-r)\Big)d\bm{Y}] 
   \mathcal{S}(\mathcal{R}(\bm{X})),
\label{eq:Hc_equality}
\end{equation}
\begin{align}
& J_i^{(c)}=H^{(c)}(\bm{X}) v_i(\bm{X})
  -\frac{\sigma^3}{2m}\int_{\mathbb{S}^2}\int_0^1
   [v_j(\bm{X}^+_{\tau\sigma\bm{\alpha}})
   -v_j(\bm{X}^+_{(\tau-1)\sigma\bm{\alpha}})]\rho(\bm{X}^+_{\tau\sigma\bm{\alpha}}) 
   \rho(\bm{X}^+_{(\tau-1)\sigma\bm{\alpha}})\notag\\
&\qquad\qquad\times  \mathcal{S}(\mathcal{R}(\bm{X}^+_{\tau\sigma\bm{\alpha}}))
   \chi_D(\bm{X}^+_{(\tau-1)\sigma\bm{\alpha}})\chi_D(\bm{X}^+_{\tau\sigma\bm{\alpha}})
   d\tau
   \alpha_i \alpha_j d\Omega(\bm{\alpha}).
\label{eq:Jc}
\end{align}%
\end{subequations}
\end{widetext}
\noindent
As far as the boundary $\partial D$ is impermeable, the integration of the second term on the right-hand side vanishes, thanks to the Gauss divergence theorem. Therefore, \eqref{eq:Hc_equality} is reduced to 
\begin{equation}
\frac{\partial {H}^{(c)}}{\partial t}+\frac{\partial J_i^{(c)}}{\partial X_i}
=-\mathcal{I}.
\label{eq:Hc0_equality}
\end{equation}
\begin{remark}
If $D$ is a unit of periodic domain or a domain surrounded by a control surface away 
from a solid boundary (or by a control surface in an infinite domain),
$\chi_D$ should be understood as unity and the spatial domain of integration occurring in the definition of $\mathcal{R}$ [see \eqref{eq:R_def}] 
should be understood to extend over $\mathbb{R}^3$. 
Then, the second term on the right-hand side of \eqref{eq:Hc_equality} vanishes.
\end{remark}

\subsection{Local H theorem\label{sec:LocalHtheorem}}

By summing up \eqref{eq:Hk_inequality} and \eqref{eq:Hc0_equality},
it holds that
\begin{subequations}\label{eq:LocalHTheorem}
\begin{align}
 \frac{\partial H}{\partial t}
+\frac{\partial J^H_i}{\partial X_i}\le 0, \label{eq:H_inequality}
\end{align}
\noindent
where
\begin{align}
H=&\langle f\ln f\rangle +H^{(c)}=H^{(k)}+H^{(c)},\\
J^H_i= &\langle \xi_i f\ln f\rangle+J_i^{(k)}+J_i^{(c)},
\end{align}
\end{subequations}
\noindent
and the equality in \eqref{eq:H_inequality} holds if and only if the condition \eqref{eq:equilibrium} is satisfied.
This is a principal part of the local H theorem for the EESM.

\begin{remark}
The flux $J^H_i$ occurring in \eqref{eq:H_inequality} is different from the flux $H_i^{(k)}+H_i^{(c)}$ defined in \cite{TT25} 
by the amount of
\begin{align}
\Delta_i\equiv &
J_i^{(k)}+J_i^{(c)}-H^{(c)}v_i \notag\\
&
-\rho v_i\int_D \frac{\rho(\bm{Y})}{m}
\theta(\sigma-|\bm{Y}-\bm{X}|)\mathcal{S}(\mathcal{R}(\bm{Y}))d\bm{Y}.\label{eq:fluxdiff}
\end{align}
%
%
\end{remark}


When the gas system is in contact with a heat bath with a uniform constant temperature $T_w$, the free energy rather than the entropy of the physical system should be a monotonic function in time.
Accordingly, it is natural to consider the local statement for the free energy. 
Following \cite{TT25}, let us introduce the function
\begin{align}
F\equiv & RT_w(\langle f\ln \frac{f}{f_w} \rangle +H^{(c)}) \notag \\
      = & RT_w H+\langle \frac12\bm{\xi}^2 f \rangle \notag \\
      = & RT_w H+ \rho(e^{(k)}+\frac12 \bm{v}^2), \label{eq:F}
\end{align}
\noindent
where 
\begin{equation}
f_{w}=\exp(-\frac{\bm{\xi}^{2}}{2RT_{w}}).\label{eq:fw}
\end{equation}
\noindent
Equation~\eqref{eq:fw} is a specific choice of a constant multiple of $f_w$ in \cite{TT25}
and gives the simplest expression for $F$.
Equation~\eqref{eq:H_inequality} is then recast as
\begin{subequations}\label{eq:LocalFreeEnergy}
\begin{align}
& \frac{\partial F}{\partial t}
+\frac{\partial J^F_i}{\partial X_i}\le 0, \label{eq:F_inequality} \\
& F= RT_w H+ \rho(e+\frac12 \bm{v}^2), \\
& J^F_i\equiv RT_w J^H_i+\rho v_i(e+\frac12\bm{v}^2)+p_{ij}v_j+q_i,\label{eq:JF}
\end{align}
\end{subequations}
\noindent
with the aid of the energy conservation \eqref{eq:eng_div},
where $e$, $p_{ij}$ and $q_i$ are the ones defined by \eqref{eq:kin+col}.
This is a secondary part of the local H theorem for the EESM.
\begin{remark}
The global quantities $\mathcal{H}^{(k)}$, $\mathcal{H}^{(c)}$, $I$, and $\mathcal{F}$ in \cite{TT25}
are nothing else than the spatial integrations of $H^{(k)}$, $H^{(c)}$, $\mathcal{I}$, and $F$ introduced in Sec.~\ref{sec:main}. Compare \eqref{H_kinetic}, \eqref{eq:H^c}, \eqref{eq:I_def}, and \eqref{eq:F} with (6), (14), (11), and (21b) in \cite{TT25}, respectively.
\end{remark}

\section{Extension to the Enskog--Vlasov equation\label{sec:EV}}

In the case of the Enskog--Vlasov equation, an external force term $F_{i}{\partial f}/{\partial\xi_{i}}$ is added on the left-hand side of \eqref{eq:2.1}, where
\begin{equation}
F_{i}=-\int_{D}\frac{\partial}{\partial X_{i}}\Phi(|\bm{Y}-\bm{X}|)\rho(\bm{Y})d\bm{Y},\label{Vlasov}
\end{equation}
\noindent
and $\Phi$ is an attractive isotropic force potential between molecules.
In the sequel, this type of external force term is referred to as the Vlasov term.

\subsection{Role of the Vlasov term in conservation laws and EoS\label{sec:VlasEoS}}

The contribution of the Vlasov term to the momentum conservation 
can be transformed into a divergence form as
\begin{align}
&\langle \xi_j F_{i}\frac{\partial f}{\partial\xi_{i}}\rangle
=\frac{\partial}{\partial X_i}p_{ij}^{(v)}, \\
&p_{ij}^{(v)}
=-\frac12 \int_{\mathbb{R}^3} \frac{r_ir_j}{|\bm{r}|}\Phi^\prime(|\bm{r}|) 
  \int_0^1 \rho(\bm{X}^+_{\lambda\bm{r}})\chi_D(\bm{X}^+_{\lambda\bm{r}})\notag\\
& \qquad\quad\times\rho(\bm{X}^+_{(\lambda-1)\bm{r}})\chi_D(\bm{X}^+_{(\lambda-1)\bm{r}})
 d\lambda d\bm{r},\label{eq:pijv}
\end{align}
\noindent
where $\Phi^\prime(x)=d\Phi(x)/dx$, which is non-negative since $\Phi$ is the attractive potential.

Similarly,
the contribution to the energy conservation can be transformed into the following form:
\begin{align}
 & \langle\frac{1}{2}\bm{\xi}^{2}F_{i}\frac{\partial f}{\partial\xi_{i}}\rangle
   = \frac{\partial}{\partial X_i}(p_{ij}^{(v)}v_j+q_i^{(v)}+\rho e^{(v)}v_i)
    +\frac{\partial}{\partial t}\rho e^{(v)}, \label{eq:VlasovEnergy}\displaybreak[0] \\
 &
   e^{(v)}=\frac12\int_D \Phi(|\bm{Y}-\bm{X}|)\rho(\bm{Y})d\bm{Y}, \label{eq:ev}\displaybreak[0]\\
 & q_i^{(v)} =(\rho e^{(v)} \delta_{ij}-p_{ij}^{(v)})v_j \notag \\
 & \qquad +\frac12\int_{\mathbb{R}^3} \int_0^1 r_i \Phi(|\bm{r}|)
           \frac{\partial\rho(\bm{X}^+_{(\lambda-1)\bm{r}})}
                {\partial t}\rho(\bm{X}^+_{\lambda\bm{r}}) \notag\\
 & \qquad\times\chi_D(\bm{X}^+_{(\lambda-1)\bm{r}})\chi_D(\bm{X}^+_{\lambda\bm{r}}) d\lambda d\bm{r}.\label{eq:qiv}
\end{align}
\noindent
The quantities  $p_{ij}^{(v)}$,  $q_i^{(v)}$, and $e^{(v)}$
can be recognized as additional contributions
to the stress tensor, the heat-flow vector, and the internal energy from the Vlasov term. 
\noindent
Hence, 
by redefining the stress tensor, heat-flow vector, and {specific} internal energy as
\begin{align}
&p_{ij}=p_{ij}^{(k)}+p_{ij}^{(c)}+p_{ij}^{(v)},\quad
q_i=q_i^{(k)}+q_i^{(c)}+q_i^{(v)},\notag\\
&e=e^{(k)}+e^{(v)},\label{eq:kin+col+vl}
\end{align}
\noindent
\eqref{eq:mom_div} and \eqref{eq:eng_div} are recovered.
Thus, together with the continuity equation \eqref{eq:continuity} that remains unchanged,
the usual form of the system of conservation equations is recovered.

In the uniform equilibrium state in the bulk,
$p_{ij}^{(v)}$ and $e^{(v)}$ are reduced to
\begin{subequations}
\begin{align}
& p_{ij}^{(v)}=-\frac{2\pi}{3}\rho^2
\int_0^\infty {x}^3 \Phi^\prime(x) dx\delta_{ij},\displaybreak[0] \\
& e^{(v)}=2\pi\rho \int_0^\infty x^2\Phi(x) dx,
\end{align}
\end{subequations}
\noindent
and thus the attractive part $p^{(v)}$ defined by
\begin{equation}
p^{(v)}=-a\rho^2,\quad
a\equiv \frac{2\pi}{3} \int_0^\infty x^3 \Phi^\prime(x) dx (>0),
\end{equation}
\noindent
is added to the right-hand side of the {EoS},
e.g., \eqref{eq:EoSvdW} and \eqref{eq:EoSCS}.
In this way, the attractive part of the {EoS}, if it exists,
is recovered by the Vlasov term.
Incidentally, as far as $x^3\Phi(x)=0$ for $x=0$ and $x\to\infty$,
$e^{(v)}$ can be rewritten as $e^{(v)}=-a\rho$
and this is consistent with the equilibrium-thermodynamic relation
$e=e_\mathrm{ideal}+\int (p-T\partial p/\partial T)/\rho^2 d\rho$,
where $e_\mathrm{ideal}=(3/2)RT(=e^{(k)})$ is the internal energy of ideal monatomic gases.

\subsection{Local H theorem\label{sec:LocalHVlasov}}

Since
\begin{equation}
\langle(1+{\ln f})F_{i}\frac{\partial f}{\partial\xi_{i}}\rangle=\langle F_{i}\frac{\partial}{\partial\xi_{i}}(f{\ln f})\rangle=0, \label{eq:Vlasovlog}
\end{equation}
\noindent
the $(1+{\ln{f}})$-moment of the Vlasov term vanishes
and thus does not contribute to \eqref{eq:27}.
Hence, the local H theorem \eqref{eq:H_inequality} remains valid as it stands.

Next consider the multiplication of the Vlasov term by $(1+\ln({f}/f_{w}))$ 
for the local statement on the free energy.
Because of \eqref{eq:Vlasovlog},
\begin{align}
\langle(1+\ln\frac{f}{f_{w}})F_{i}\frac{\partial f}{\partial\xi_{i}}\rangle 
&=-\langle({\ln{f_{w}}})F_{i}\frac{\partial f}{\partial\xi_{i}}\rangle \notag\\
&=\frac{1}{RT_w}\langle\frac{1}{2}\bm{\xi}^{2}F_{i}\frac{\partial f}{\partial\xi_{i}}\rangle,
\label{eq:externalF}
\end{align}
\noindent
and the Vlasov term has a new contribution, 
which is rewritten by using \eqref{eq:VlasovEnergy} as
\begin{align}
&\langle(1+\ln\frac{f}{f_{w}})F_{i}\frac{\partial f}{\partial\xi_{i}}\rangle
\Big(
=\frac{1}{RT_w}\langle\frac{1}{2}\bm{\xi}^{2}F_{i}\frac{\partial f}{\partial\xi_{i}}\rangle\Big)
  \notag \\
&=\frac{1}{RT_w}
    \Big( \frac{\partial}{\partial t}\rho e^{(v)}
    +\frac{\partial}{\partial X_i}[p_{ij}^{(v)}v_j+q_i^{(v)}+\rho e^{(v)}v_i]
    \Big).
    \label{eq:externalF2}
\end{align}
\noindent
Therefore, the local statement on the free energy is modified as
\begin{subequations}\label{eq:tildeFinequality}
\begin{equation}
\frac{\partial \widetilde{F}}{\partial t}+\frac{\partial J^{\widetilde{F}}_i}{\partial X_i}\le 0,
\end{equation}
\noindent
holds, where
\begin{align}
 \widetilde{{F}}
\equiv & {F}+\rho e^{(v)},\\
J^{\widetilde{F}}_i
=&J^{{F}}_i + p_{ij}^{(v)}v_j+q_i^{(v)}+\rho e^{(v)}v_i.
\end{align}
\end{subequations}

\begin{remark}
The difference between \eqref{eq:LocalFreeEnergy} and \eqref{eq:tildeFinequality} 
is confined in the stress tensor, heat-flow vector, and specific internal energy.
Therefore, by switching from \eqref{eq:kin+col} to \eqref{eq:kin+col+vl}, 
\eqref{eq:LocalFreeEnergy} covers the modified version \eqref{eq:tildeFinequality}. 
\end{remark}

\section{Conclusion\label{sec:Conclusion}}

In the present paper,
the local H theorem has been shown to hold for EESM and EVESM proposed in \cite{TT25}.
This result refines the corresponding global statement in \cite{TT25}
in the sense that it is pointwise and thus is a more detailed statement.
Thanks to the refinement,
the local H function that has a direct link to the local entropy production is indeed identified 
together with the newly defined fluxes. 
Hence the present results are expected to serve as the basis of the reciprocity arguments 
for EESM and EVESM.

\appendix
\section{Some standard operations for the collision term
and their consequences \label{app:shift}}
We summarize some standard operations
that are used in the transformations of the collision term.
There are three types of operation that
are standard in the case of the Boltzmann equation as well:
\begin{description}
\item [{(I)}] to exchange the letters $\bm{\xi}$ and $\bm{\xi}_{*}$; 
\item [{(II)}] to reverse the direction of $\bm{\alpha}$;
\item [{(III)}] to change the integration variables from $(\bm{\xi},\bm{\xi}_{*},\bm{\alpha})$
to $(\bm{\xi}^{\prime},\bm{\xi}_{*}^{\prime},\bm{\alpha})$ and then to
change the letters $(\bm{\xi}^{\prime},\bm{\xi}_{*}^{\prime})$ to
$(\bm{\xi},\bm{\xi}_{*})$.
\end{description}

First, by (III) and (II),
\begin{equation}
\int{\varphi}(\bm{X}) J^{G}(f)d\bm{\xi}=\int{\varphi}^{\prime}(\bm{X})J^{L}(f)d\bm{\xi},\label{eq:4.1}
\end{equation}
\noindent
holds for any ${\varphi}(\bm{X},\bm{\xi})$.
Hence, we have
\begin{align}
 & \langle\varphi(\bm{X})J(f)\rangle\nonumber \\
=& \frac{\sigma^{2}}{m}
   \int[{\varphi}^\prime(\bm{X})-{\varphi}(\bm{X})]g(\bm{X}_{\sigma\bm{\alpha}}^{-},\bm{X})f_{*}(\bm{X}_{\sigma\bm{\alpha}}^{-})f(\bm{X}) \notag\\ &\qquad\times V_{\alpha}\theta(V_{\alpha})d\Omega(\bm{\alpha})d\bm{\xi}_{*}d\bm{\xi}.\label{eq:mom_J}
\end{align}

\subsection{Transformation of $\langle J(f)\ln f \rangle$}

The substitution of $\varphi=\ln f$ in \eqref{eq:mom_J} yields
\begin{align}
 & \langle J(f)\ln f\rangle\nonumber \\
=& \frac{\sigma^{2}}{m}
   \int[\ln f^\prime(\bm{X})-\ln f(\bm{X})]g(\bm{X}_{\sigma\bm{\alpha}}^{-},\bm{X})f_{*}(\bm{X}_{\sigma\bm{\alpha}}^{-})f(\bm{X}) \notag\\ &\qquad\times V_{\alpha}\theta(V_{\alpha})d\Omega(\bm{\alpha})d\bm{\xi}_{*}d\bm{\xi}\nonumber \\
=& \frac{\sigma^{2}}{2m}
   \int[\ln \Big(\frac{f^\prime(\bm{X})f^\prime_{*}(\bm{X}^-_{\sigma\bm{\alpha}})}
                      {f(\bm{X})f_{*}(\bm{X}^-_{\sigma\bm{\alpha}})}
            \Big)
       +\ln \Big(\frac{f^\prime(\bm{X})f_{*}(\bm{X}^-_{\sigma\bm{\alpha}})}
                      {f(\bm{X})f^\prime_{*}(\bm{X}^-_{\sigma\bm{\alpha}})}
            \Big)
        ]
    \notag\\
 &\qquad\times g(\bm{X}_{\sigma\bm{\alpha}}^{-},\bm{X})f_{*}(\bm{X}_{\sigma\bm{\alpha}}^{-})f(\bm{X}) 
      V_{\alpha}\theta(V_{\alpha})d\Omega(\bm{\alpha})d\bm{\xi}_{*}d\bm{\xi}\nonumber \\
=& \frac{\sigma^{2}}{2m}
   \int\ln \Big(\frac{f^\prime(\bm{X})f^\prime_{*}(\bm{X}^-_{\sigma\bm{\alpha}})}
                      {f(\bm{X})f_{*}(\bm{X}^-_{\sigma\bm{\alpha}})}
            \Big)
    \notag\\
 &\qquad\times g(\bm{X}_{\sigma\bm{\alpha}}^{-},\bm{X})f_{*}(\bm{X}_{\sigma\bm{\alpha}}^{-})f(\bm{X}) 
      V_{\alpha}\theta(V_{\alpha})d\Omega(\bm{\alpha})d\bm{\xi}_{*}d\bm{\xi}\nonumber \\
 &-\frac{\sigma^{2}}{2m}
   \int[\ln \frac{f^\prime_{*}(\bm{X}^-_{\sigma\bm{\alpha}})}{f_{*}(\bm{X}^-_{\sigma\bm{\alpha}})}
       -\ln \frac{f^\prime(\bm{X})}{f(\bm{X})}]
    \notag\\
 &\qquad\times g(\bm{X}_{\sigma\bm{\alpha}}^{-},\bm{X})f_{*}(\bm{X}_{\sigma\bm{\alpha}}^{-})f(\bm{X}) 
      V_{\alpha}\theta(V_{\alpha})d\Omega(\bm{\alpha})d\bm{\xi}_{*}d\bm{\xi}.
\end{align}
\noindent
But, the second term on the right-hand side is further transformed as
\begin{align}
 &-\frac{\sigma^{2}}{2m}
   \int[\ln \frac{f^\prime_{*}(\bm{X}^-_{\sigma\bm{\alpha}})}{f_{*}(\bm{X}^-_{\sigma\bm{\alpha}})}
       -\ln \frac{f^\prime(\bm{X})}{f(\bm{X})}]
    \notag\\
 &\qquad\times g(\bm{X}_{\sigma\bm{\alpha}}^{-},\bm{X})f_{*}(\bm{X}_{\sigma\bm{\alpha}}^{-})f(\bm{X}) 
      V_{\alpha}\theta(V_{\alpha})d\Omega(\bm{\alpha})d\bm{\xi}_{*}d\bm{\xi} \notag
      \displaybreak[0] \\
=&-\frac{\sigma^{2}}{2m}
   \int[\ln \frac{f^\prime(\bm{X}^+_{\sigma\bm{\alpha}})}{f(\bm{X}^+_{\sigma\bm{\alpha}})}
   g(\bm{X}_{\sigma\bm{\alpha}}^{+},\bm{X})f(\bm{X}_{\sigma\bm{\alpha}}^{+})f_{*}(\bm{X}) \notag \\
& \qquad -\ln \frac{f^\prime(\bm{X})}{f(\bm{X})}
       g(\bm{X}_{\sigma\bm{\alpha}}^{-},\bm{X})f_{*}(\bm{X}_{\sigma\bm{\alpha}}^{-})f(\bm{X}) ] \notag\\
& \qquad \times V_{\alpha}\theta(V_{\alpha})d\Omega(\bm{\alpha})d\bm{\xi}_{*}d\bm{\xi} 
\notag \displaybreak[0]\\
=&-\frac{\partial}{\partial X_j}\Big(
   \frac{\sigma^{3}}{2m}
   \int\int_0^1 \alpha_j g(\bm{X}_{\tau\sigma\bm{\alpha}}^{+},\bm{X}^{+}_{(\tau-1)\sigma\bm{\alpha}}) \notag\\
&  \qquad \times [\ln \frac{f^\prime(\bm{X}^+_{\tau\sigma\bm{\alpha}})}{f(\bm{X}^+_{\tau\sigma\bm{\alpha}})}
   f(\bm{X}_{\tau\sigma\bm{\alpha}}^{+})f_{*}(\bm{X}^+_{(\tau-1)\sigma\bm{\alpha}})]d\tau
    \notag\\
 &\qquad\times V_{\alpha}\theta(V_{\alpha})d\Omega(\bm{\alpha})d\bm{\xi}_{*}d\bm{\xi} 
   \Big),   
\label{eq:mom_Jlog2}
\end{align}
where the first transformation is the result of applying operations (I) and (II) to the first term.
Consequently, the following relation has been obtained:
\begin{subequations}\label{eq:JlnJ}
\begin{align}
 & \langle J(f)\ln f\rangle\nonumber \\
=& \frac{\sigma^{2}}{2m}
   \int g(\bm{X}_{\sigma\bm{\alpha}}^{-},\bm{X}) 
   \ln \Big(\frac{f^\prime(\bm{X})f^\prime_{*}(\bm{X}^-_{\sigma\bm{\alpha}})}
                      {f(\bm{X})f_{*}(\bm{X}^-_{\sigma\bm{\alpha}})}
            \Big)f_{*}(\bm{X}_{\sigma\bm{\alpha}}^{-})f(\bm{X})
    \displaybreak[0]\notag\\
 &\qquad\times  
      V_{\alpha}\theta(V_{\alpha})d\Omega(\bm{\alpha})d\bm{\xi}_{*}d\bm{\xi}
      -\frac{\partial J_j^{(k)}}{\partial X_j},   
\label{eq:mom_Jlog}
\end{align}
\noindent
where
\begin{align}
  J_j^{(k)}
=&\frac{\sigma^{3}}{2m}
   \int\int_0^1 \alpha_j g(\bm{X}_{\tau\sigma\bm{\alpha}}^{+},\bm{X}^{+}_{(\tau-1)\sigma\bm{\alpha}}) \displaybreak[0]\notag\\
&  \quad \times [\ln \frac{f^\prime(\bm{X}^+_{\tau\sigma\bm{\alpha}})}{f(\bm{X}^+_{\tau\sigma\bm{\alpha}})}
   f(\bm{X}_{\tau\sigma\bm{\alpha}}^{+})f_{*}(\bm{X}^+_{(\tau-1)\sigma\bm{\alpha}})]d\tau
    \notag\\
 &\quad\times V_{\alpha}\theta(V_{\alpha})d\Omega(\bm{\alpha})d\bm{\xi}_{*}d\bm{\xi}. 
 \label{eq:Jk}
\end{align}
\end{subequations}
\noindent   
Equation~\eqref{eq:mom_Jlog} is the form of \eqref{eq:Jlnf}.

\subsection{Transformation of $\mathcal{I}$ \label{app:ColEnt}}

The transformation of $\mathcal{I}$ from \eqref{eq:I_def} to \eqref{eq:I_final}
has been done by operation (III) followed by (II):
\begin{align}
\mathcal{I} 
= & -\frac{\sigma^{2}}{2m}\int g(\bm{X},\bm{X}_{\sigma\bm{\alpha}}^{-})f_{*}^{\prime}(\bm{X}_{\sigma\bm{\alpha}}^{-})f^{\prime}(\bm{X})  \notag\\ 
 &\qquad\qquad\times V_{\alpha}^{\prime}\theta(-V_{\alpha}^{\prime}) d\Omega(\bm{\alpha})d\bm{\xi}d\bm{\xi}_{*}\nonumber  \\
 & -\frac{\sigma^{2}}{2m}\int g(\bm{X},\bm{X}_{\sigma\bm{\alpha}}^{-})f(\bm{X})f_{*}(\bm{X}_{\sigma\bm{\alpha}}^{-}) \notag\\
 &\qquad\qquad\times V_{\alpha}\theta(V_{\alpha}) d\Omega(\bm{\alpha})d\bm{\xi}d\bm{\xi}_{*}\displaybreak[0]\nonumber \\
= & -\frac{\sigma^{2}}{2m}\int g(\bm{X},\bm{X}_{\sigma\bm{\alpha}}^{-})f(\bm{X})f_{*}(\bm{X}_{\sigma\bm{\alpha}}^{-}) \notag\\
 &\qquad\qquad\times V_{\alpha} d\Omega(\bm{\alpha})d\bm{\xi}d\bm{\xi}_{*}\displaybreak[0]\nonumber \\
= &\frac{\sigma^{2}}{2m}\int g(\bm{X},\bm{X}_{\sigma\bm{\alpha}}^{-})
 \rho(\bm{X})\rho(\bm{X}_{\sigma\bm{\alpha}}^{-})\notag\\
&\qquad
\times[ v_j(\bm{X})-v_j(\bm{X}_{\sigma\bm{\alpha}}^{-})]\alpha_j 
 d\Omega(\bm{\alpha})\displaybreak[0]\nonumber \\
= &\frac{\sigma^{2}}{2m}\int g(\bm{X},\bm{X}_{\sigma\bm{\alpha}}^{+})
 \rho(\bm{X})\rho(\bm{X}_{\sigma\bm{\alpha}}^{+})\notag\\
&\qquad
\times[ v_j(\bm{X}_{\sigma\bm{\alpha}}^{+})-v_j(\bm{X})]\alpha_j 
 d\Omega(\bm{\alpha})
,\label{eq:5.4}
\end{align}
\noindent
where $V_{\alpha}^{\prime}\equiv(\bm{\xi}_{*}^{\prime}-\bm{\xi}^{\prime})\cdot\bm{\alpha}=-V_{\alpha}$ has been used
at the beginning of the above transformation.

\section{Connection with the global statement\label{app:mono}}

In \cite{TT25}, the global H theorem has been discussed,
especially for three typical cases: the domain $D$ is three dimensional and 
is (i) periodic, (ii) surrounded by the specular reflection boundary,
and (iii) surrounded by the impermeable surface of a heat bath
with a uniform constant temperature $T_w$. 
We will show that the monotonic decrease of the global function
\begin{equation}
\mathcal{H}\equiv \int_D H d\bm{X},
\end{equation}
is recovered in cases (i) and (ii),
and that the monotonic decrease of another global function
\begin{equation}
\mathcal{F}\equiv \int_D F d\bm{X},
\end{equation}
is recovered in case (iii).

In order to see the monotonicity of $\mathcal{H}$, 
let us focus on the difference of flux
$\Delta_i$ given in \eqref{eq:fluxdiff}.

First, the divergence of the last term on the right-hand side of \eqref{eq:fluxdiff}
vanishes after the spatial integration over $D$ and using the Gauss divergence theorem and no mass flux condition $v_in_i=0$ on the impermeable boundary, namely in cases (ii) and (iii).
This vanishing remains true in case (i), 
since the surface integral on $\partial D$ cancels out by periodic condition.

Next, consider the divergence of $\mathcal{J}_i^{(c)}\equiv J_i^{(c)}-H^{(c)}v_i$.
It holds from \eqref{eq:Jcori} and \eqref{eq:Jc} that 
\begin{widetext}
\begin{align}
 & \int_{D} \frac{\partial}{\partial X_i}\mathcal{J}_i^{(c)}d\bm{X} \notag\\
=&-\int_{\mathbb{R}^3}
   \frac{\sigma^2}{2m}\int_{\mathbb{S}^2}
   \Big(
   [v_i(\bm{X}^+_{\sigma\bm{\alpha}})
   -v_i(\bm{X})]\rho(\bm{X}^+_{\sigma\bm{\alpha}}) 
   \rho(\bm{X})\mathcal{S}(\mathcal{R}(\bm{X}^+_{\sigma\bm{\alpha}}))
   \chi_D(\bm{X})\chi_D(\bm{X}^+_{\sigma\bm{\alpha}})
   \notag \displaybreak[0]\\
& -[v_i(\bm{X})-v_i(\bm{X}^-_{\sigma\bm{\alpha}})]
   \rho(\bm{X}^-_{\sigma\bm{\alpha}}) 
   \rho(\bm{X})\mathcal{S}(\mathcal{R}(\bm{X}))
   \chi_D(\bm{X})\chi_D(\bm{X}^-_{\sigma\bm{\alpha}})
   \Big)\alpha_i
   d\Omega(\bm{\alpha}) d\bm{X} \notag \displaybreak[0]\\
=&-\frac{\sigma^2}{2m}\int_{\mathbb{S}^2}
   \int_{\mathbb{R}^3}
   \Big(
   [v_i(\bm{X}^+_{\sigma\bm{\alpha}})-v_i(\bm{X})]
   \rho(\bm{X}^+_{\sigma\bm{\alpha}}) 
   \rho(\bm{X})\mathcal{S}(\mathcal{R}(\bm{X}^+_{\sigma\bm{\alpha}}))
   \chi_D(\bm{X})\chi_D(\bm{X}^+_{\sigma\bm{\alpha}})
   \Big)\alpha_i  d\bm{X} d\Omega(\bm{\alpha})
   \notag \displaybreak[0]\\
& + \frac{\sigma^2}{2m}\int_{\mathbb{S}^2}
   \int_{\mathbb{R}^3}
   [v_i(\bm{Z}^+_{\sigma\bm{\alpha}})-v_i(\bm{Z})]
   \rho(\bm{Z})\rho(\bm{Z}^+_{\sigma\bm{\alpha}}) 
   \mathcal{S}(\mathcal{R}(\bm{Z}^+_{\sigma\bm{\alpha}}))
   \chi_D(\bm{Z}^+_{\sigma\bm{\alpha}})\chi_D(\bm{Z})
   \Big)\alpha_i d\bm{Z} d\Omega(\bm{\alpha}),
\end{align}
\end{widetext}
where the last transformation is done by setting $\bm{Z}=\bm{X}^-_{\sigma\bm{\alpha}}$.
Obviously, the two terms on the most right-hand side cancel out each other.
The divergence of $\mathcal{J}_i^{(c)}$ thus vanishes once integrated over the domain $D$.

Finally, consider the divergence of $J_i^{(k)}$.
As is seen from \eqref{eq:Jk} and \eqref{eq:mom_Jlog2}, it holds that
\begin{align}
 &\int_D \frac{\partial}{\partial X_i}J_i^{(k)}d\bm{X}\displaybreak[0] \notag \\
=&\int_D
   \frac{\sigma^{2}}{2m}
   \int[\ln \frac{f^\prime(\bm{X}^+_{\sigma\bm{\alpha}})}{f(\bm{X}^+_{\sigma\bm{\alpha}})}
   g(\bm{X}_{\sigma\bm{\alpha}}^{+},\bm{X})f(\bm{X}_{\sigma\bm{\alpha}}^{+})f_{*}(\bm{X}) \displaybreak[0]\notag \\
& \qquad -\ln \frac{f^\prime(\bm{X})}{f(\bm{X})}
       g(\bm{X}_{\sigma\bm{\alpha}}^{-},\bm{X})f_{*}(\bm{X}_{\sigma\bm{\alpha}}^{-})f(\bm{X}) ] \displaybreak[0]\notag\\
& \qquad \times V_{\alpha}\theta(V_{\alpha})d\Omega(\bm{\alpha})d\bm{\xi}_{*}d\bm{\xi}
  d\bm{X} \displaybreak[0]\notag\\
=&\frac{\sigma^{2}}{2m}
   \int\int_{\mathbb{R}^3}[ \ln \frac{f^\prime(\bm{X}^+_{\sigma\bm{\alpha}})}{f(\bm{X}^+_{\sigma\bm{\alpha}})}
   g(\bm{X}_{\sigma\bm{\alpha}}^{+},\bm{X})f(\bm{X}_{\sigma\bm{\alpha}}^{+})f_{*}(\bm{X}) \displaybreak[0]\notag \\
& \qquad -\ln \frac{f^\prime(\bm{X})}{f(\bm{X})}
       g(\bm{X}_{\sigma\bm{\alpha}}^{-},\bm{X})f_{*}(\bm{X}_{\sigma\bm{\alpha}}^{-})f(\bm{X}) ] d\bm{X} \displaybreak[0]\notag\\
& \qquad \times V_{\alpha}\theta(V_{\alpha})d\Omega(\bm{\alpha})d\bm{\xi}_{*}d\bm{\xi} \displaybreak[0]\notag\\
=&\frac{\sigma^{2}}{2m}
   \int[\int_{\mathbb{R}^3} \ln \frac{f^\prime(\bm{X}^+_{\sigma\bm{\alpha}})}{f(\bm{X}^+_{\sigma\bm{\alpha}})}
   g(\bm{X}_{\sigma\bm{\alpha}}^{+},\bm{X})f(\bm{X}_{\sigma\bm{\alpha}}^{+})f_{*}(\bm{X})d\bm{X} \displaybreak[0]\notag \\
& \qquad -\int_{\mathbb{R}^3}\ln \frac{f^\prime(\bm{Z}^+_{\sigma\bm{\alpha}})}{f(\bm{Z}^+_{\sigma\bm{\alpha}})}
       g(\bm{Z},\bm{Z}^+_{\sigma\bm{\alpha}})f_{*}(\bm{Z})f(\bm{Z}^+_{\sigma\bm{\alpha}})  d\bm{Z} ]\displaybreak[0]\notag\\
& \qquad \times V_{\alpha}\theta(V_{\alpha})d\Omega(\bm{\alpha})d\bm{\xi}_{*}d\bm{\xi},
\end{align}
\noindent
similar to the case of $\mathcal{J}_i^{(c)}$.
Obviously again, the last two terms cancel out each other.
Thus the divergence of $J_i^{(k)}$ vanishes as well once integrated over the domain $D$.

Therefore the difference of the flux from the one in \cite{TT25} vanishes 
after the integration over the domain $D$. 
In this way, the global statement for the monotonic decrease of $\mathcal{H}$ in \cite{TT25} is recovered
for cases (i) and (ii).
Since the difference of $J_i^F$ from $J^H_i$ originates from $\langle\xi_i f\ln f_w\rangle$,
the difference can be handled by the Darrozes--Guiraud inequality~\cite{DG66,C88,S07} for case (iii). 
Hence, the monotonic decrease of $\mathcal{F}$ for case (iii)
results along the lines of the argument in \cite{TT25}.

In closing, let us examine the influence of the Vlasov term.
Firstly, since the local H theorem is not affected by the Vlasov term,
$\mathcal{H}$ decreases monotonically in time even in the presence of the Vlasov term for cases (i) and (ii), as proved in \cite{TT25}.
Secondly, because of \eqref{eq:externalF2} and  
\begin{equation}
 \frac{1}{RT_w}\langle\frac{1}{2}\bm{\xi}^{2}F_{i}\frac{\partial f}{\partial\xi_{i}}\rangle
=-\frac{\rho v_i F_i}{RT_w},
\end{equation}
\noindent
it holds that
\begin{align}
&-\int_D {\rho v_i F_i}d\bm{X} \notag
\\
=&\int_D [\frac{\partial}{\partial t}\rho e^{(v)}
+ \frac{\partial}{\partial X_i} (p_{ij}^{(v)}v_j+q_i^{(v)}+\rho e^{(v)}v_i)] d\bm{X}.
\end{align}
\noindent
In the meantime, according to \cite{TT25} [see (E4) in its Appendix~E],
it also holds that
\begin{equation}
-\int_D \rho v_i F_i d\bm{X}
=\frac{d}{dt}\int_D \rho e^{(v)} d\bm{X}.
\end{equation}
\noindent
Hence, it holds that
\begin{equation}
 \int_D \frac{\partial}{\partial X_i} (p_{ij}^{(v)}v_j+q_i^{(v)}+\rho e^{(v)}v_i) d\bm{X}
=0.
\end{equation}
\noindent
Therefore, the divergence of the additional flux $(p_{ij}^{(v)}v_j+q_i^{(v)}+\rho e^{(v)}v_i)$ 
in $J^{\tilde{F}}_{i}$ (against $J^F_i$) vanishes, once integrated over the domain $D$. The monotonic decrease of 
\begin{equation}
\widetilde{\mathcal{F}}=\int_D \widetilde{F}d\bm{X},
\end{equation}
\noindent
is thus recovered for case~(iii)~\cite{errata}.

\acknowledgments
The present work is partially supported by the JSPS Grant-in-Aid for JSPS Fellows (No.~24KJ1450) for the first author and by the Kyoto University Foundation and the JSPS Grant-in-Aid for Scientific Research(B) (No.~26K00870) for the second author.


\begin{thebibliography}{99}
\bibitem{errata}
The description about $\widetilde{H}$ at the beginning of page~10 of \cite{TT25} is not appropriate.
The part ``, in place of ... or'' in the 2nd to 5th lines of the left column on the same page should be deleted. Accordingly, the last line of Appendix E should read ``... Appendix D is extended to that of $\mathcal{H}$ and $\widetilde{\mathcal{F}}$ ...''
\end{thebibliography}
\end{document}